\documentclass[a4paper,12pt,oneside]{article}
\usepackage{amsmath,amsfonts,amssymb,amsmath,latexsym,amsthm,enumerate}

\newtheorem{thm}{Theorem}
\newtheorem{cor}[thm]{Corollary}
\theoremstyle{definition}

\theoremstyle{plain}

\begin{document}
\title {Umbral Calculus and Frobenius-Euler Polynomials}
\author{by \\Dae San Kim and Taekyun Kim}\date{}\maketitle
\begin{abstract}
\noindent In this paper, we study some properties of umbral calculus related to Appell sequence. From those properties, we derive new and interesting identities of Frobenius-Euler polynomials.
\end{abstract}
\section{Introduction}
Let $\mathbf{C}$ be the complex number field. For $\lambda\in\mathbf{C}$ with $\lambda\neq 1$, the Frobenius-Euler polynomials are defined by the generating function to be
\begin{equation}\label{eq:1}
\frac{1-\lambda}{e^{t}-\lambda}e^{xt}=e^{H(x\vert\lambda)t}=\sum_{n=0}^{\infty}H_{n}(x\vert\lambda)\frac{t^{n}}{n!},\quad (\text{see}\,\, \lbrack 7\!\!-\!\!11\rbrack)\,,
\end{equation}
with the usual convention about replacing $H^{n}(x\vert\lambda)$ by $H_{n}(x\vert\lambda)$.\\
In the special case, $x=0, H_{n}(0\vert\lambda)=H_{n}(\lambda)$ are called the $n$-th Frobenius-Euler numbers.
By (\ref{eq:1}), we get
\begin{equation}\label{eq:2}
H_{n}(x\vert\lambda)=\sum_{l=0}^{n}\left(\begin{array}{c}n\\l\end{array}\right)H_{n-l}(\lambda)x^{l}=(H(\lambda)+x)^{n}, \quad(\text{see}\,\,\lbrack1,2,3,4\rbrack),
\end{equation}
with the usual convention about replacing $H^{n}(\lambda)$ by $H_{n}(\lambda)$.\\
Thus, from (\ref{eq:1}) and (\ref{eq:2}), we note that
\begin{equation*}
(H(\lambda)+1)^{n}-\lambda H_{n}(\lambda)=(1-\lambda)\delta_{0,n}\,\,,
\end{equation*}
where $\delta_{n,k}$ is the kronecker symbol (see\,\,\lbrack6,7\rbrack).\\
For $r\in\mathbf{Z}_{+}$, the Frobenius-Euler polynomials of order $r$ are defined by the generating function to be
\begin{align}\label{eq:3}
(\frac{1-\lambda}{e^{t}-\lambda})e^{xt}&=\underbrace{(\frac{1-\lambda}{e^{t}-\lambda})\times\cdots\times(\frac{1-\lambda}{e^{t}-\lambda})e^{xt}}_{r-times}\\
                                      &=\sum_{n=0}^{\infty}H_{n}^{(r)}(x\vert\lambda)\frac{t^{n}}{n!}.\nonumber
\end{align}
In the special case, $x=0$, $H_{n}^{(r)}(0\vert\lambda)=H_{n}^{(r)}(\lambda)$ are called the $n$-th Frobenius-Euler numbers of order $r$ (see\,\,\lbrack6,7\rbrack).\\
From (\ref{eq:3}), we can derive the following equation:
\begin{equation}\label{eq:4}
H_{n}^{(r)}(x\vert\lambda)=\sum_{l=0}^{n}\left(\begin{array}{c}n\\l\end{array}\right)H_{n-l}^{(r)}(\lambda)x^{l},
\end{equation}
and
\begin{equation}\label{eq:5}
H_{n}^{(r)}(\lambda)=\sum_{l_{1}+\cdots+l_{r}=n}\left(\begin{array}{c}n\\l_{1},\cdots,l_{r}\end{array}\right)H_{l_{1}}(\lambda)\cdots H_{l_{r}}(\lambda).
\end{equation}
By (\ref{eq:4}) and (\ref{eq:5}), we see that $H_{n}^{(r)}(x\vert\lambda)$ is a monic polynomial of degree $n$ with coefficients in $\mathbf{Q}(\lambda)$.\\
Let $\mathbb{P}$ be the algebra of polynomials in the single variable $x$ over $\mathbf{C}$ and let $\mathbb{P}^{*}$ be the vector space of all linear functionals on $\mathbb{P}$. As is known, $\langle L\vert p(x)\rangle$ denotes the action of the linear functional $L$ on a polynomial $p(x)$ and we remind that the addition and scalar multiplication on $\mathbb{P}^{*}$ are respectively defined by
\begin{equation*}
\langle L+M\vert p(x)\rangle=\langle L\vert p(x)\rangle +\langle M\vert p(x)\rangle, \langle cL\vert p(x) \rangle = c\langle L\vert p(x)\rangle,
\end{equation*}
where $c$ is a complex constant (see\,\,\cite{DS,Rom}).\\
Let $\mathbf{F}$ denote the algebra of formal power series:
\begin{equation}\label{eq:6}
\mathbf{F}=\lbrace f(t)=\sum_{k=0}^{\infty}\frac{a_{k}}{k!}t^{k}\vert a_{k}\in\mathbf{C}\rbrace,\quad (\text{see}\,\,\cite{DS,Rom}).
\end{equation}
The formal power series define a linear functional on $\mathbb{P}$ by setting
\begin{equation}\label{eq:7}
\langle f(t)\vert x^{n}\rangle = a_{n},\,\, \text{for all}\,\, n\geq0.
\end{equation}
Indeed, by (\ref{eq:6}) and (\ref{eq:7}), we get
\begin{equation}\label{eq:8}
\langle t^{k}\vert x^{n} \rangle = n!\delta_{n,k}\quad (n,k\geq 0),\,\,\,\,(\text{see}\,\, \cite{DS,Rom}).
\end{equation}
This kind of algebra is called an umbral algebra.\\
The order $O(f(t))$ of a nonzero power series $f(t)$ is the smallest integer $k$ for which the coefficient of $t^{k}$ does not vanish. A series $f(t)$ for which $O(f(t))=1$ is said to be an invertible series (see $\lbrack5,8\rbrack$). For $f(t), g(t)\in\mathbf{F}$ and $p(x)\in\mathbb{P}$, we have
\begin{equation}\label{eq:9}
\langle f(t)g(t)\vert p(x)\rangle = \langle f(t)\vert g(t)p(x)\rangle=\langle g(t)\vert f(t)p(x)\rangle,\,\, (\text{see}\,\,\cite{DS}).
\end{equation}
One should keep in mind that each $f(t)\in\mathbf{F}$ plays three roles in the umbral calculus : a formal power series, a linear functional and a linear operator. To illustrate this, let $p(x)\in\mathbb{P}$ and $f(t)=e^{yt}\in\mathbf{F}$. As a linear functional, $e^{yt}$ satisfies $\langle e^{yt}\vert p(x) \rangle=p(y)$.
As a linear operator, $e^{yt}$ satisies $e^{yt}p(x)=p(x+y)$ (see $\lbrack5\rbrack$). Let $s_{n}(x)$ denote a polynomial in $x$ with degree $n$. Let us assume that $f(t)$ is a delta series and $g(t)$ is an invertible series. Then there exists a unique sequence $s_{n}(x)$ of polynomials such that $\langle g(t)f(t)^{k}\vert s_{n}(x)\rangle=n!\delta_{n,k}$ for all $n,k\geq 0$ (see\,\, \cite{DS,Rom}). This sequence $s_{n}(x)$ is called the Sheffer sequence for $(g(t), f(t))$ which is denoted by $s_{n}(x)\sim(g(t), f(t))$. If $s_{n}(x)\sim(1, f(t))$, then $s_{n}(x)$ is called the associated sequence for $f(t)$. If $s_{n}(x)\sim(g(t), t)$, then $s_{n}(x)$ is called the Appell sequence.\\
Let $s_{n}(x)\sim(g(t),f(t))$. Then we see that
\begin{align}
\label{eq:10} h(t)&=\sum_{k=0}^{\infty}\frac{<h(t)\vert s_{k}(x)>}{k!}g(t)f(t)^{k},\quad h(t)\in\mathbf{F}\,,\\
\label{eq:11} p(x)&=\sum_{k=0}^{\infty}\frac{<g(t)f(t)^{k}\vert p(x)>}{k!}s_{k}(x),\quad p(x)\in\mathbb{P}\,,\\
\label{eq:12} f(t)s_{n}(x)&=ns_{n-1}(x),\quad <f(t)\vert p(\alpha x)>=<f(\alpha t\vert p(x)>,
\end{align}
and
\begin{equation}\label{eq:13}
\frac{1}{g(\bar{f}(t))}e^{y\bar{f}(t)}=\sum_{k=0}^{\infty}\frac{s_{k}(y)}{k!}t^{k},\quad\text{for all}\,\,y\in\mathbf{C},
\end{equation}
where $\bar{f}(t)$ is the compositional inverse of $f(t)$ (see\,\,\cite{Rom}).
In this paper, we study some properties of umbral calculus related to Appell sequence. For those properties, we derive new and interesting of Frobenius-Euler polynomials.
 \section{Frobenius-Euler polynomials and Umbral Calculus.}
 By (\ref{eq:3}) and (\ref{eq:13}), we see that
 \begin{equation}\label{eq:14}
 H_{n}^{(r)}(x\vert\lambda)\sim((\frac{e^{t}-\lambda}{1-\lambda})^{r}, t)\,.
 \end{equation}
 Thus, by (\ref{eq:14}), we get
 \begin{equation}\label{eq:15}
 <(\frac{e^{t}-\lambda}{1-\lambda})^{r}t^{k}\vert H_{n}^{(r)}(x\vert\lambda)>=n!\delta_{n,k}\,\,.
 \end{equation}
 Let
 \begin{equation*}
 \mathbb{P}_{n}(\lambda)=\lbrace p(x)\in\mathbf{Q}(\lambda)\lbrack x\rbrack\vert\,\, \text{deg}\,\, p(x)\leq n\rbrace\,.
 \end{equation*}
 Then it is an $(n+1)$-dimensional vector space over $\mathbf{Q}(\lambda)$.\\
 So we see that $\lbrace H_{0}^{(r)}(x\vert\lambda), H_{1}^{(r)}(x\vert\lambda), \cdots, H_{n}^{(r)}(x\vert\lambda)\rbrace$ is a basis for $\mathbb{P}_{n}(\lambda)$. For $p(x)\in\mathbb{P}_{n}(\lambda)$, let
 \begin{equation}\label{eq:16}
 p(x)=\sum_{k=0}^{n}C_{k}H_{k}^{(r)}(x\vert\lambda),\quad (n\geq 0).
 \end{equation}
 Then, by (\ref{eq:14}), (\ref{eq:15}) and (\ref{eq:16}), we get
 \begin{align}\label{eq:17}
 <(\frac{e^{t}-\lambda}{1-\lambda})^{r}t^{k}\vert p(x)>&=\sum_{l=0}^{n}C_{l}<(\frac{e^{t}-\lambda}{1-\lambda})^{r}t^{k}\vert H_{l}^{(r)}(x\vert\lambda)>\\
                                                       &=\sum_{l=0}^{n}C_{l}l!\delta_{l,k}=k!C_{k}\,. \nonumber
 \end{align}
 From (\ref{eq:17}), we have
 \begin{align}\label{eq:18}
 C_{k}&=\frac{1}{k!}<(\frac{e^{t}-\lambda}{1-\lambda})^{r}t^{k}\vert p(x)>=\frac{1}{k!}<(\frac{e^{t}-\lambda}{1-\lambda})^{r}\vert D^{k}p(x)>\\
      &=\frac{1}{k!(1-\lambda)^{r}}\sum_{j=0}^{r}\left(\begin{array}{c}r\\j\end{array}\right)(-\lambda)^{r-j}<e^{jt}\vert D^{k}p(x)>\nonumber\\
      &=\frac{1}{k!(1-\lambda)^{r}}\sum_{j=0}^{r}\left(\begin{array}{c}r\\j\end{array}\right)(-\lambda)^{r-j}<t^{0}\vert e^{jt}D^{k}p(x)>\nonumber\\
      &=\frac{1}{k!(1-\lambda)^{r}}\sum_{j=0}^{r}\left(\begin{array}{c}r\\j\end{array}\right)(-\lambda)^{r-j}<t^{0}\vert D^{k}p(x+j)>\,.\nonumber
 \end{align}
Therefore, by (\ref{eq:16}) and (\ref{eq:18}), we obtain the following theorem.
\begin{thm}\label{thm1}
For $p(x)\in\mathbb{P}_{n}(\lambda)$\,, let
\begin{equation*}
p(x)=\sum_{k=0}^{n}C_{k}H_{k}^{(r)}(x)\,.
\end{equation*}
Then we have
\begin{equation*}
C_{k}=\frac{1}{k!(1-\lambda)^{r}}\sum_{j=0}^{r}\left(\begin{array}{c}r\\j\end{array}\right)(-\lambda)^{r-j}D^{k}p(j)\,\,,
\end{equation*}
where $Dp(x)=\frac{dp(x)}{dx}$.
\end{thm}
\noindent From Theorem \ref{thm1}, we note that
\begin{equation*}
p(x)=\frac{1}{(1-\lambda)^{r}}\sum_{k=0}^{n}\lbrace\sum_{j=0}^{r}\frac{1}{k!}\left(\begin{array}{c}r\\j\end{array}\right)(-\lambda)^{r-j}D^{k}p(j)\rbrace H_{k}^{(r)}(x\vert\lambda)\,.
\end{equation*}
Let us consider the operator $\tilde{\triangle}_{\lambda}$ with $\tilde{\triangle}_{\lambda}f(x)=f(x+1)-\lambda f(x)$ and let $J_{\lambda}=\frac{1}{1-\lambda}\tilde{\triangle}_{\lambda}$\,. Then we have
\begin{equation}\label{eq:19}
J_{\lambda}(f)(x)=\frac{1}{1-\lambda}\lbrace f(x+1)-\lambda f(x)\rbrace.
\end{equation}
Thus, by (\ref{eq:19}), we get
\begin{equation}\label{eq:20}
J_{\lambda}(H_{n}^{(r)}(x\vert\lambda)=\frac{1}{1-\lambda}\lbrace H_{n}^{(r)}(x+1\vert\lambda)-\lambda H_{n}^{(r)}(x\vert\lambda)\rbrace.
\end{equation}
From (\ref{eq:3}), we can derive
\begin{align}\label{eq:21}
&\sum_{n=0}^{\infty}\lbrace H_{n}^{(r)}(x+1\vert\lambda)-\lambda H_{n}^{(r)}(x\vert\lambda)\rbrace\frac{t^{n}}{n!}=(\frac{1-\lambda}{e^{t}-\lambda})^{r}e^{(x+1)t}-\lambda(\frac{1-\lambda}{e^{t}-\lambda})e^{xt}\\
&=(\frac{1-\lambda}{e^{t}-\lambda})^{r}e^{xt}(e^{t}-\lambda)=(1-\lambda)(\frac{1-\lambda}{e^{t}-\lambda})^{r-1}e^{xt}=(1-\lambda)\sum_{n=0}^{\infty}H_{n}^{(r-1)}(x\vert\lambda)\frac{t^{n}}{n!}.\nonumber
\end{align}
By (\ref{eq:20}) and (\ref{eq:21}), we get
\begin{equation}\label{eq:22}
J_{\lambda}(H_{n}^{(r)}(x\vert\lambda))=H_{n}^{(r-1)}(x\vert\lambda).
\end{equation}
From (\ref{eq:22}), we have
\begin{equation*}
J_{\lambda}^{r}(H_{n}^{(r)}(x\vert\lambda))=J_{\lambda}^{r-1}(H_{n}^{(r-1)}(x\vert\lambda))=\cdots=H_{n}^{(0)}(x\vert\lambda)=x^{n},
\end{equation*}
and
\begin{equation}\label{eq:23}
J_{\lambda}^{r}(x^{n})=J_{\lambda}^{r}H_{n}^{(0)}(x\vert\lambda)=H_{n}^{(-r)}(x\vert\lambda)=J_{\lambda}^{2r}H_{n}^{(r)}(x\vert\lambda)\,.
\end{equation}
For $s\in\mathbf{Z}_{+}$, from (\ref{eq:22}), we have
\begin{equation}\label{eq:24}
J_{\lambda}^{s}(H_{n}^{(r)}(x\vert\lambda))=H_{n}^{(r-s)}(x\vert\lambda)\,.
\end{equation}
On the other hand, by (\ref{eq:13}), (\ref{eq:14}) and (\ref{eq:22}),
\begin{align}\label{eq:25}
J_{\lambda}^{s}(H_{n}^{(r)}(x\vert\lambda))&=(\frac{e^{t}-\lambda}{1-\lambda})^{s}(H_{n}^{(r)}(x\vert\lambda))\\
                                           &=\frac{1}{(1-\lambda)^{s}}((1-\lambda)+\sum_{k=1}^{\infty}\frac{t^{k}}{k!})^{s}(H_{n}^{(r)}(x\vert\lambda)).\nonumber
\end{align}
Thus, by\,\,(\ref{eq:25}), we get
\begin{align}
&\label{eq:26}J_{\lambda}^{s}(H_{n}^{(r)}(x\vert\lambda)=\sum_{m=0}^{s}\frac{\left(\begin{array}{c}s\\m\end{array}\right)}{(1-\lambda)^{m}}\sum_{l=m}^{\infty}(\sum_{k_1+\cdots+k_m=l\atop k_j\geq1}\frac{1}{k_1!\cdots k_m!})t^{l}(H_{n}^{(r)}(x\vert\lambda))\\
&=\sum_{m=0}^{s}\frac{\left(\begin{array}{c}s\\m\end{array}\right)}{(1-\lambda)^{m}}\sum_{l=m}^{\infty}\frac{1}{l!}(\sum_{k_1+\cdots+k_m=l\atop k_j\geq1}\left(\begin{array}{c}l\\k_{1},\cdots,k_{m}\end{array}\right)D^{l}(H_{n}^{(r)}(x\vert\lambda))\nonumber\\
&=\sum_{m=0}^{min\lbrace s,n \rbrace} \frac{\left(\begin{array}{c}s\\m\end{array}\right)}{(1-\lambda)^{m}}\sum_{l=m}^{n}\left(\begin{array}{c}n\\l\end{array}\right)\sum_{k_1+\cdots+k_m=l\atop k_j\geq1}\left(\begin{array}{c}l\\k_{1},\cdots,k_{m}\end{array}\right)H_{n-l}^{(r)}(x\vert\lambda)\nonumber\\
&=\sum_{l=0}^{min\lbrace s,n \rbrace}\lbrace\left(\begin{array}{c}n\\l\end{array}\right)\sum_{m=0}^{l}\frac{\left(\begin{array}{c}s\\m\end{array}\right)}{(1-\lambda)^{m}}\sum_{k_1+\cdots+k_m=l\atop k_j\geq1}\left(\begin{array}{c}l\\k_{1},\cdots,k_{m}\end{array}\right)\rbrace H_{n-l}^{(r)}(x\vert\lambda)\nonumber\\
&+\sum_{l=min\lbrace s, n\rbrace+1}^{n}\lbrace\left(\begin{array}{c}n\\l\end{array}\right)\sum_{m=0}^{min\lbrace s,n\rbrace}\frac{\left(\begin{array}{c}s\\m\end{array}\right)}{(1-\lambda)^{m}}\sum_{k_1+\cdots+k_m=l\atop k_j\geq1}\left(\begin{array}{c}l\\k_{1},\cdots,k_{m}\end{array}\right)\rbrace H_{n-l}^{(r)}(x\vert\lambda)\nonumber
\end{align}
Therefore, by (\ref{eq:24}) and (\ref{eq:26}), we obtain the following theorem.
\begin{thm}\label{theorem2}
For any $r, s\geq 0$, we have
\begin{align*}
&H_{n}^{(r-s)}(x\vert\lambda)=\sum_{l=0}^{min\lbrace s,n\rbrace}\lbrace\left(\!\!\begin{array}{c}n\\l\end{array}\!\!\right)\sum_{m=0}^{l}\frac{\left(\begin{array}{c}s\\m\end{array}\right)}{(1-\lambda)^{m}}\sum_{k_1+\cdots+k_m=l\atop k_j\geq1}\left(\!\!\begin{array}{c}l\\k_{1},\cdots,k_{m}\end{array}\!\!\right)\rbrace H_{n-l}^{(r)}(x\vert\lambda)\\
&+\sum_{l=min\lbrace s,n\rbrace+1}^{n}\lbrace\left(\begin{array}{c}n\\l\end{array}\right)\sum_{m=0}^{min\lbrace s,n\rbrace}\frac{\left(\begin{array}{c}s\\m\end{array}\right)}{(1-\lambda)^{m}}\sum_{k_1+\cdots+k_m=l\atop k_j\geq1}\left(\begin{array}{c}l\\k_{1},\cdots,k_{m}\end{array}\right)\rbrace H_{n-l}^{(r)}(x\vert\lambda).
\end{align*}
\end{thm}
\noindent Let us take $s=r-1 (r\geq1)$ in Theorem \ref{theorem2}. Then we obtain the following corollary.
\begin{cor}\label{cor3}
For $n\geq0, r\geq1,$ we have
\begin{align*}
&H_{n}(x\vert\lambda)\\
&=\sum_{l=0}^{min\lbrace r-1, n\rbrace}\lbrace\left(\begin{array}{c}n\\l\end{array}\right)\sum_{m=0}^{l}\frac{\left(\begin{array}{c}r-1\\m\end{array}\right)}{(1-\lambda)^{m}}\sum_{k_1+\cdots+k_m=l\atop k_j\geq1}\left(\begin{array}{c}l\\k_{1},\cdots,k_{m}\end{array}\right)\rbrace H_{n-l}^{(r)}(x\vert\lambda)\\
&+\!\!\sum_{l=min\lbrace r-1,n\rbrace+1}^{n}\lbrace\left(\!\!\begin{array}{c}n\\l\end{array}\!\!\right)\sum_{m=0}^{min\lbrace r-1,n\rbrace}\frac{\left(\!\!\begin{array}{c}r-1\\m\end{array}\!\!\right)}{(1-\lambda)^{m}}\sum_{k_1+\cdots+k_m=l\atop k_j\geq1}\left(\!\!\begin{array}{c}l\\k_{1},\cdots,k_{m}\end{array}\!\!\right)\rbrace H_{n-l}^{(r)}(x\vert\lambda).
\end{align*}
\end{cor}
\noindent Let us take $s=r (r\geq1)$ in Theorem \ref{theorem2}. Then we obtain the following corollary.
\begin{cor}\label{cor4}
For $n\geq0, r\geq1$, we have
\begin{align*}
x^{n}&=\sum_{l=0}^{min\lbrace r, n\rbrace}\lbrace\left(\!\begin{array}{c}n\\l\end{array}\!\right)\sum_{m=0}^{l}\frac{\left(\begin{array}{c}r\\m\end{array}\right)}{(1-\lambda)^{m}}\sum_{k_1+\cdots+k_m=l\atop k_j\geq1}\left(\!\!\begin{array}{c}l\\k_{1},\cdots,k_{m}\end{array}\!\!\right)\rbrace H_{n-l}^{(r)}(x\vert\lambda)\\
&+\!\sum_{l=min\lbrace r,n\rbrace+1}^{n}\lbrace\left(\begin{array}{c}n\\l\end{array}\right)\sum_{m=0}^{min\lbrace r,n\rbrace}\frac{\left(\begin{array}{c}r\\m\end{array}\right)}{(1-\lambda)^{m}}\sum_{k_1+\cdots+k_m=l\atop k_j\geq1}\left(\!\!\begin{array}{c}l\\k_{1},\cdots,k_{m}\end{array}\!\!\right)\rbrace H_{n-l}^{(r)}(x\vert\lambda).
\end{align*}
\end{cor}
\noindent Now, we define the analogue of Stirling numbers of the second kind as follows:
\begin{equation}\label{eq:27}
S_{\lambda}(n,k)=\frac{1}{k!}\sum_{j=0}^{k}\left(\begin{array}{c}k\\j\end{array}\right)(-\lambda)^{k-j}j^{n},\,\,\,\,(n,k\geq0)\,.
\end{equation}
Note that $S_{1}(n,k)=S(n,k)$ is the stirling number of the second kind.\\
From the definition of $\tilde{\triangle}_{\lambda}$, we have
\begin{equation}\label{eq:28}
\tilde{\triangle}_{\lambda}^{n}f(0)=\sum_{k=0}^{n}\left(\begin{array}{c}n\\k\end{array}\right)(-\lambda)^{n-k}f(k)
\end{equation}
By (\ref{eq:27}) and (\ref{eq:28}), we get
\begin{equation}\label{eq:29}
S_{\lambda}(n,k)=\frac{1}{k!}\tilde{\triangle}_{\lambda}^{k}0^{n},\,\,(n,k\geq0).
\end{equation}
Let us take $s=2r$. Then we have
\begin{align}
&\label{eq:30}J_{\lambda}^{r}x^{n}=H_{n}^{(-r)}(x\vert\lambda)\\
                                 &=\sum_{l=0}^{min\lbrace2r,n\rbrace}\lbrace \left(\begin{array}{c}n\\l\end{array}\right)\sum_{m=0}^{l}\frac{\left(\begin{array}{c}2r\\m\end{array}\right)}{(1-\lambda)^{m}}\sum_{k_1+\cdots+k_{m}=l\atop k_j\geq1}\left(\begin{array}{c}l\\k_1,\cdots,k_m\end{array}\right)\rbrace H_{n-l}^{(r)}(x\vert\lambda)\nonumber\\
                                 &+\sum_{l=min\lbrace2r,n\rbrace+1}^{n}\!\lbrace\left(\begin{array}{c}n\\l\end{array}\right)\!\sum_{m=0}^{min\lbrace2r,n\rbrace}\frac{\left(\begin{array}{c}2r\\m\end{array}\right)}{(1-\lambda)^{m}}\sum_{k_1+\cdots+k_m=l\atop k_j\geq1}\left(\begin{array}{c}l\\k_1,\cdots,k_m\end{array}\right)\rbrace H_{n-l}^{(r)}(x\vert\lambda),\nonumber
\end{align}
and
\begin{equation}\label{eq:31}
J_{\lambda}^{r}x^{n}=(\frac{1}{1-\lambda}\tilde{\triangle}_{\lambda})^{r}(x^{n})=\frac{1}{(1-\lambda)^{r}}\sum_{j=0}^{r}\left(\begin{array}{c}r\\j\end{array}\right)(-\lambda)^{r-j}(x+j)^{n}.
\end{equation}
By (\ref{eq:30}) and (\ref{eq:31}), we get
\begin{align}\label{eq:32}
&\frac{1}{(1-\lambda)^{r}}\sum_{j=0}^{r}\left(\begin{array}{c}r\\j\end{array}\right)(-\lambda)^{r-j}(x+j)^{n}=\frac{1}{(1-\lambda)^{r}}\tilde{\triangle}_{\lambda}^{r}x^{n}\\
                                                                                                         &=\sum_{l=0}^{min\lbrace2r,n\rbrace}\lbrace\left(\begin{array}{c}n\\l\end{array}\right)\sum_{m=0}^{l}\frac{\left(\begin{array}{c}2r\\m\end{array}\right)}{(1-\lambda)^{m}}\sum_{k_1+\cdots+k_{m}=l\atop k_j\geq1}\left(\begin{array}{c}l\\k_1,\cdots,k_m\end{array}\right)\rbrace H_{n-l}^{(r)}(x\vert\lambda)\nonumber\\
                                                                                                         &+\sum_{l=min\lbrace2r,n\rbrace+1}^{n}\lbrace\!\left(\begin{array}{c}n\\l\end{array}\right)\!\sum_{m=0}^{min\lbrace2r,n\rbrace}\frac{\left(\begin{array}{c}2r\\m\end{array}\right)}{(1-\lambda)^{m}}\sum_{k_1+\cdots+k_m=l\atop k_j\geq1}\left(\begin{array}{c}l\\k_1,\cdots,k_m\end{array}\right)\rbrace H_{n-l}^{(r)}(x\vert\lambda).\nonumber
\end{align}
Let us take $x=0$ in (\ref{eq:32}). Then we obtain the following theorem.
\begin{thm}\label{theorem5}
\begin{align*}
&\frac{r!}{(1-\lambda)^{r}}S_{\lambda}(n,r)=\frac{r!}{(1-\lambda)^{r}}\frac{\tilde{\triangle}_{\lambda}^{r}0^{n}}{r!}\\
&\sum_{l=0}^{min\lbrace 2r,n\rbrace}\lbrace \left(\begin{array}{c}n\\l\end{array}\right)\sum_{m=0}^{l}\frac{\left(\begin{array}{c}2r\\m\end{array}\right)}{(1-\lambda)^{m}}\sum_{k_1+\cdots+k_{m}=l\atop k_j\geq1}\left(\begin{array}{c}l\\k_1,\cdots,k_m\end{array}\right)\rbrace H_{n-l}^{(r)}(\lambda)\\
                                 &+\sum_{l=min\lbrace 2r, n\rbrace+1}^{n}\!\!\left(\begin{array}{c}n\\l\end{array}\right)\!\!\sum_{m=0}^{min\lbrace2r,n\rbrace}\frac{\left(\begin{array}{c}2r\\m\end{array}\right)}{(1-\lambda)^{m}}\sum_{k_1+\cdots+k_m=l\atop k_j\geq1}\left(\begin{array}{c}l\\k_1,\cdots,k_m\end{array}\right)\rbrace H_{n-l}^{(r)}(\lambda)\\
                                 &=\sum_{m=0}^{min\lbrace r,n\rbrace}\frac{\left(\begin{array}{c}r\\m\end{array}\right)}{(1-\lambda)^{m}}\sum_{k_1+\cdots+k_m=n\atop k_j\geq1}\left(\begin{array}{c}n\\k_1,\cdots,k_m\end{array}\right).
\end{align*}
\end{thm}
Let us consider $s=2r-1$ in the identity of Theorem \ref{theorem2}. Then we have
\begin{align}\label{eq:33}
&J_{\lambda}^{r-1}x^{n}=H_{n}^{-(r-1)}(x\vert\lambda)\\
                      &=\sum_{l=0}^{min\lbrace 2r-1,n\rbrace}\lbrace\left(\begin{array}{c}n\\l\end{array}\right)\sum_{m=0}^{l}\frac{\left(\begin{array}{c}2r-1\\m\end{array}\right)}{(1-\lambda)^{m}}\sum_{k_1+\cdots+k_{m}=l\atop k_j\geq1}\left(\begin{array}{c}l\\k_1,\cdots,k_m\end{array}\right)\rbrace H_{n-l}^{(r)}(x\vert\lambda)\nonumber\\
                      &+\!\!\sum_{l=min\lbrace 2r-1, n\rbrace+1}^{n}\!\lbrace\!\left(\begin{array}{c}n\\l\end{array}\right)\!\sum_{m=0}^{min\lbrace2r-1,n\rbrace}\!\frac{\left(\begin{array}{c}2r-1\\m\end{array}\right)}{(1-\lambda)^{m}}\!\sum_{k_1+\cdots+k_m=l\atop k_j\geq1}\!\left(\begin{array}{c}l\\k_1,\cdots,k_m\end{array}\right)\!\rbrace\!
                      H_{n-l}^{(r)}(x\vert\lambda)\nonumber\\
                      &=\frac{1}{(1-\lambda)^{r-1-j}}\sum_{j=0}^{r-1}\left(\begin{array}{c}r-1\\j\end{array}\right)(-\lambda)^{r-1-j}(x+j)^{n}=\frac{1}{(1-\lambda)^{r-1}}\tilde{\triangle}_{\lambda}^{r-1}x^{n}.\nonumber
\end{align}
Let us take $x=0$ in (\ref{eq:33}). Then we obtain the following theorem.
\begin{thm}\label{theorem6}
For $n\geq0$ and $r\geq1$, we have
\begin{align*}
&\frac{(r-1)!}{(1-\lambda)^{r-1}}S_{\lambda}(n,r-1)=\frac{(r-1)!}{(1-\lambda)^{r-1}}\frac{\tilde{\triangle}_{\lambda}^{r-1}0^{n}}{(r-1)!}\\
                                                  &=\sum_{l=0}^{min\lbrace 2r-1,n\rbrace}\lbrace\left(\begin{array}{c}n\\l\end{array}\right)\sum_{m=0}^{l}\frac{\left(\begin{array}{c}2r-1\\m\end{array}\right)}{(1-\lambda)^{m}}\sum_{k_1+\cdots+k_{m}=l\atop k_j\geq1}\left(\begin{array}{c}l\\k_1,\cdots,k_m\end{array}\right)\rbrace H_{n-l}^{(r)}(\lambda)\nonumber\\
                                                  &+\!\!\sum_{l=min\lbrace 2r-1, n\rbrace+1}^{n}\lbrace\!\left(\begin{array}{c}n\\l\end{array}\right)\!\sum_{m=0}^{min\lbrace2r-1,n\rbrace}\!\!\!\frac{\left(\!\!\begin{array}{c}2r-1\\m\end{array}\!\!\right)}{(1-\lambda)^{m}}\!\!\sum_{k_1+\cdots+k_m=l\atop k_j\geq1}\left(\!\!\begin{array}{c}l\\k_1,\cdots,k_m\end{array}\!\!\right)\rbrace H_{n-l}^{(r)}(\lambda)\nonumber
\end{align*}
\end{thm}
\noindent$\large{\mathbf{Remark}}$. Note that
\begin{align*}
&\frac{(r-1)!}{(1-\lambda)^{r-1}}S_{\lambda}(n,r-1)\\
&=\sum_{l=0}^{min\lbrace r,n\rbrace}\lbrace\left(\begin{array}{c}n\\l\end{array}\right)\sum_{m=0}^{l}\frac{\left(\begin{array}{c}r\\m\end{array}\right)}{(1-\lambda)^{m}}\sum_{k_1+\cdots+k_{m}=l\atop k_j\geq1}\left(\begin{array}{c}l\\k_1,\cdots,k_m\end{array}\right)\rbrace H_{n-l}(\lambda)\nonumber\\
                                                 &+\sum_{l=min\lbrace r, n\rbrace+1}^{n}\left(\begin{array}{c}n\\l\end{array}\right)\sum_{m=0}^{min\lbrace r,n\rbrace}\frac{\left(\begin{array}{c}r\\m\end{array}\right)}{(1-\lambda)^{m}}\sum_{k_1+\cdots+k_m=l\atop k_j\geq1}\left(\begin{array}{c}l\\k_1,\cdots,k_m\end{array}\right)\rbrace H_{n-l}(\lambda)\nonumber
\end{align*}

\author{Department of Mathematics, Sogang University, Seoul 121-742, Republic of Korea
\\e-mail: dskim@sogang.ac.kr}\\
\\
\author{Department of Mathematics, Kwangwoon University, Seoul 139-701, Republic of Korea
\\e-mail: tkkim@kw.ac.kr}
\end{document}